\documentclass[12pt,a4paper]{article}

	\usepackage{fancyhdr,amsmath, amssymb, amsfonts,dsfont, subfigure,amsthm,tikz,color,epsfig,graphicx,soul,mathtools}

    \usepackage{IEEEtrantools}

	\usetikzlibrary{calc}

\newlength{\guillotine}
	\newtheorem{thm}{Theorem}[section]
	
	\newtheorem{lem}[thm]{Lemma}
	
	\newtheorem{prop}[thm]{Proposition}

	\newtheorem{defn}[thm]{Definition}
	\newtheorem{example}[thm]{Example}

	\theoremstyle{remark}
	\newtheorem{rem}[thm]{Remark}

	\newcommand\RR{\mathbb R}
	\newcommand\NN{\mathbb N}
	
	\newcommand\CC{\mathbb C}

	\newcommand\cala{\mathcal A}
	\newcommand\cals{\mathcal S}
	\newcommand\calr{\mathcal R}
	
	\newcommand\calg{\mathcal G}
	
	\newcommand\cali{\mathcal I}

    \newcommand\alf\alpha
	\newcommand\de\delta
	\newcommand\eps\varepsilon
	\newcommand\la\lambda
	\newcommand\De\Delta
	\newcommand\La\Lambda
	\newcommand\Om\Omega
	\newcommand\kap\kappa
    \newcommand\sig\sigma

  	\newcommand\calo{\mathcal O}
	
	\makeatletter
	\newcommand{\ostar}{\mathbin{\mathpalette\make@circled\star}}
	\newcommand{\make@circled}[2]{%
	  \ooalign{$\m@th#1\smallbigcirc{#1}$\cr\hidewidth$\m@th#1#2$\hidewidth\cr}%
	}
	\newcommand{\smallbigcirc}[1]{%
	  \vcenter{\hbox{\scalebox{0.77778}{$\m@th#1\bigcirc$}}}%
	}
	\makeatother

	\newcommand\emp\emptyset

    \newcommand\un\underline
    
    \newcommand\ol\overline

	\DeclareMathOperator\area{area}
	
	\DeclareMathOperator\per{per}
	\DeclareMathOperator\In{in}
	\DeclareMathOperator\vol{vol}
	\DeclareMathOperator\Int{int}
	\DeclareMathOperator\nab\nabla

\title{
A formula for the upper box-counting dimension of self-projective sets
}
\author{Benedict Sewell\thanks{Supported generously by the Alfr\'ed R\'enyi Young Researcher Fund.}}
\date{\today}

\begin{document}

\maketitle

\begin{abstract}
	We prove a packing exponent--style formula for the upper box-counting dimension of attractors of certain projective iterated function systems. This partially affirms a conjecture of De Leo \cite{de-leo conjecture original}, and gives that the box-counting dimension of the Rauzy gasket $\calr$, $\dim_B(\calr)$ exists and satisfies
		$$
				\dim_B(\calr)
			=
				\dim_H(\calr)
			\in
				\big[
					1.6196,
					1.7415
				\big]
			.
		$$
	%
\end{abstract}

\section{Introduction}

The Rauzy gasket is an interesting set.
Discovered independently no fewer than three times, it is an important subset of parameter space for numerous problems in dynamics and topology, in fields as diverse as crystalline conductivity \cite{maltsev-novikov}, triangular tiling billiards \cite{hubert-paris romaskevich}, and asymptotics on Markov--Hurwitz varieties \cite{gamburd-magee-ronan}. Outside of these contexts, it is interesting in its own right as a fractal object, specifically for the challenge it poses: concrete bounds on its fractal dimensions (specifically Hausdorff and box-counting) have been hard to come by, despite a lot of recent attention.
The gasket is a prototypical example of a self-projective attractor. The fractal geometric theory of such sets has so far only been undertaken in one dimension \cite{barnsley-vince,christodolou-jurga,solomyak-takahashi}, with the recent exception of \cite{chinese}.

As suggested by the title, our main result is an explicit formula for the upper box-counting dimension of self-projective attractors with simplicial holes, which we prove using tools from classical analytic number theory. This allows us to derive an lower bound supporting a general conjecture of De Leo \cite{de-leo conjecture original}, and we can also give an upper bound in the more general case of convex holes. Concerning the Rauzy gasket $\calr$, the main result applies, together with a corresponding formula of \cite{chinese}, to show that its box-counting dimension $\dim_B(\calr)$ exists and equals its Hausdorff dimension $\dim_H(\calr)$, and to give the strongest numerical lower bound for this (Hausdorff) dimension to date.

Before we state our results, we introduce the Rauzy gasket, illustrating some of the difficulty concerned by comparing it to two other well-known gaskets, and listing what is known about its dimension.

\subsection{Good, bad and ugly gaskets}

The three gaskets we touch upon in this section are mutually homeomorphic \cite{arnoux-starosta} but extremely different in the fractal geometric sense. Such differences can even be seen visually in the figures given below.

The word gasket\footnote{We much prefer the French term, \textit{baderne}---thanks to P. Arnoux for the introduction.} in each case refers to the attractor $\calg$ of three distance-decreasing maps $f_k$ of some triangular subset of the plane, i.e., $\calg$ is the unique non-empty compact set such that
	$
			\calg
		=
			f_1(\calg)
		\cup
			f_2(\calg)
		\cup
			f_3(\calg).
	$

Certain properties of the $f_j$, specifically whether they are weakly or strongly contracting, and how distorting they are, are reflected in the relative ease of determining, or estimating, the fractal dimensions of the corresponding $\calg$, as we shall see.

\subsubsection{The Sierpi\'nski gasket}

\begin{figure}[ht]
	\centering
	\includegraphics[width=0.5\linewidth]{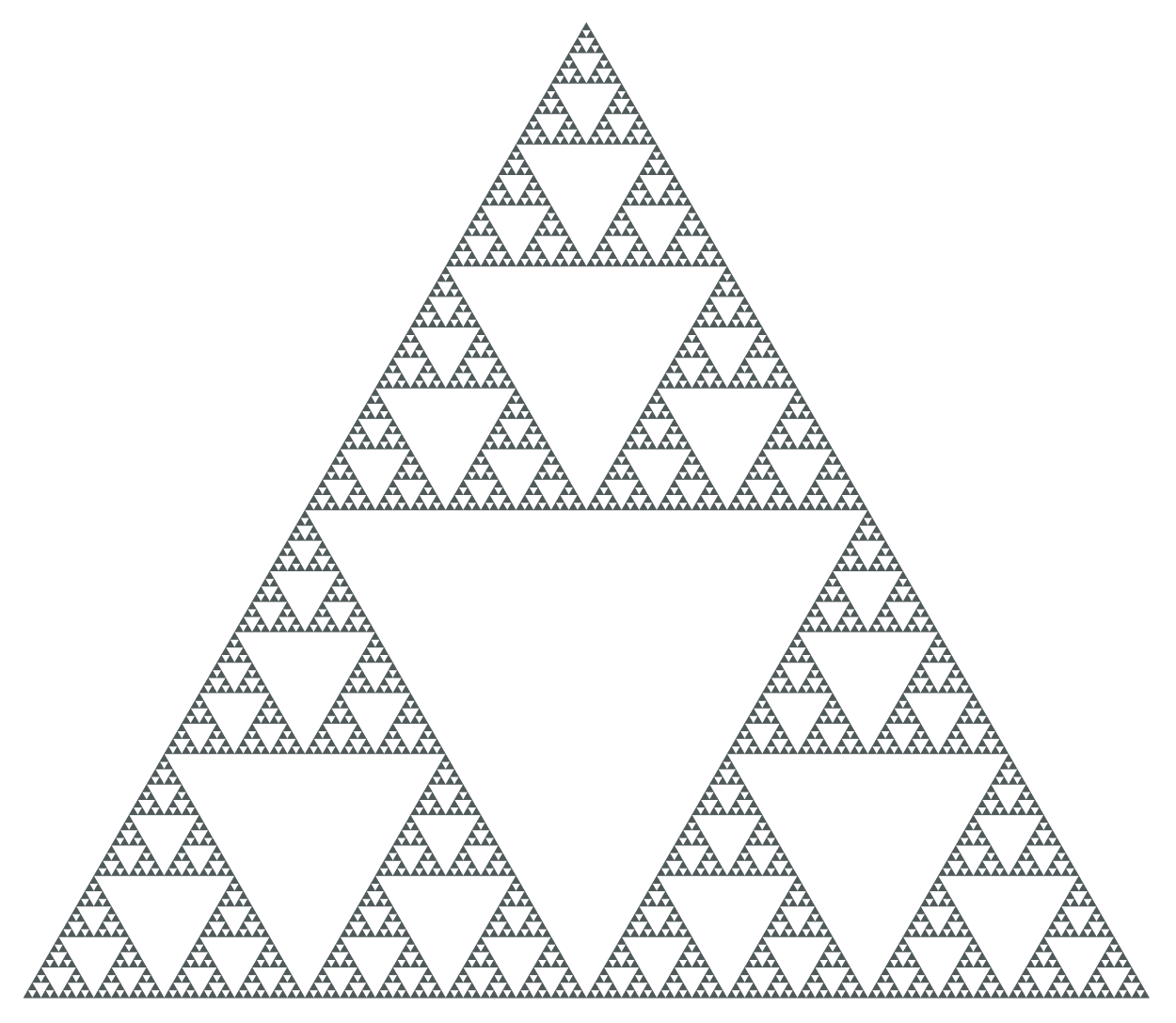}
	\caption[The Sierpi\'nski gasket $\cals$]{The Sierpi\'nski gasket, $\cals$.}
	\label{fig:sierpinski}
\end{figure}

The Sierpi\'nski gasket, $\cals$, introduced by Wac\l{}aw Sierpi\'nski in 1915 \cite{sierpinski} and shown in Figure \ref{fig:sierpinski}, is the example seen in first courses in fractal geometry.

Up to linear deformation, it is the limit set of three contracting similarities $x \mapsto \frac 12(x + e_k)$, where $(e_k)_{k=1}^3 \subset \RR^2$ are fixed, non-colinear points. The rigidity and strong contractiveness of these maps makes the Hausdorff and box-counting dimensions easy to explicitly calculate:
	$$
			\dim_H(S)
		=
			\dim_B(S)
		=
			\frac{\log(3)}{\log(2)}
		=
			1.5849\ldots.
	$$
Hence it is not particularly interesting from the fractal geometric viewpoint. We however note some interesting connections to, e.g., cellular automata and Pascal's triangle, which can be found on Wikipedia.

\subsubsection{The Apollonian gasket}

\begin{figure}[ht]
	\centering
	\includegraphics[width=0.9\linewidth]{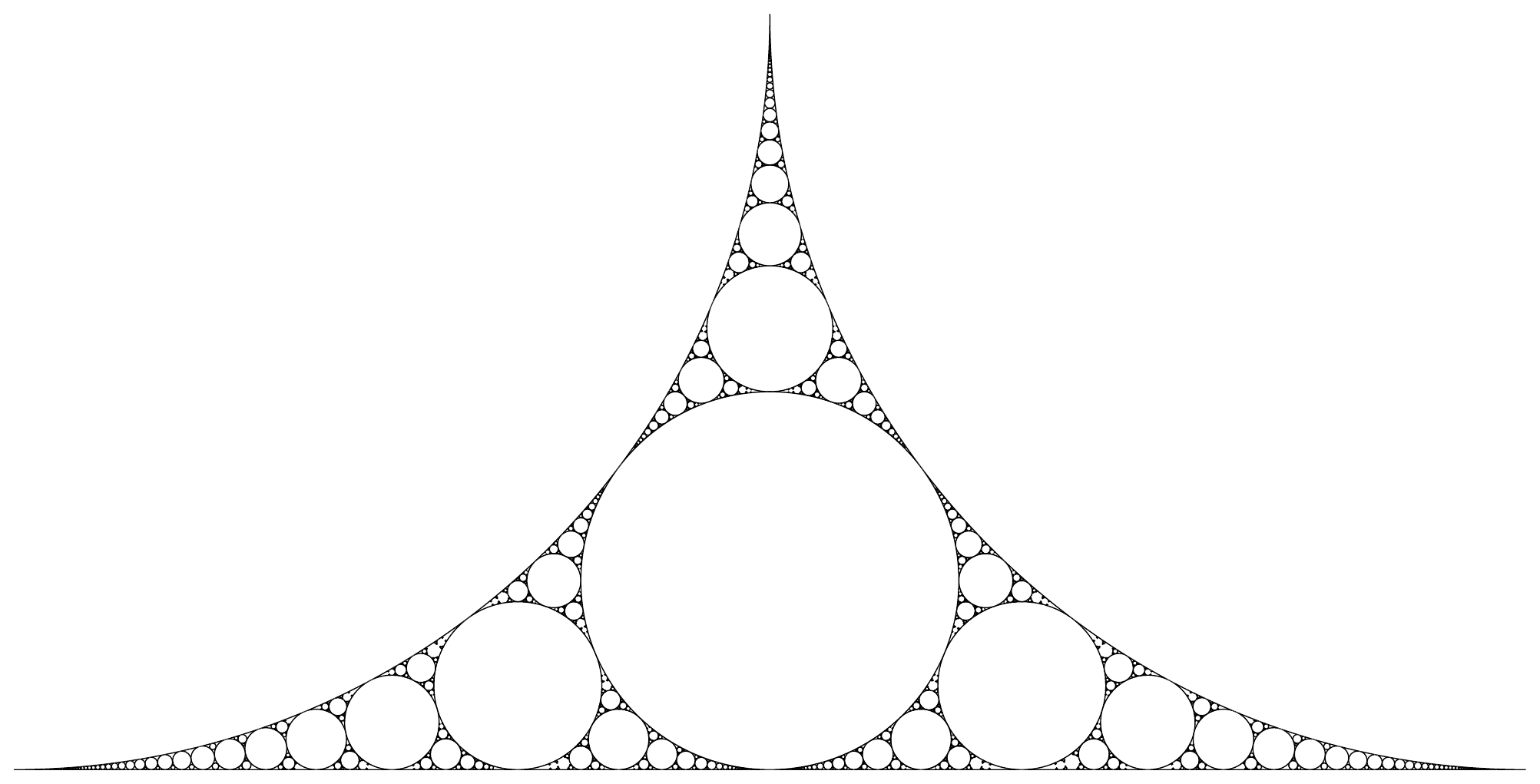}
	\caption[The Apollonian gasket $\cala$]{The Apollonian gasket, $\cala$}
	\label{figr-apollonian}
\end{figure}

The Apollonian gasket $\cala$ has perhaps the longest and richest history of any fractal, dating back to a theorem of Apollonius of Perga (ca.~200BC), which even features a formula of Decartes. We defer the reader to \cite{pollicott} for a readable historical account.

Up to a conformal change of coordinates, $\cala$ is the attractor of three M\"obius maps:
	$$
		z \mapsto \frac {z-1}{z+3},
	\qquad
		 z \mapsto \frac {z+1}{3-z},
	\qquad{\text{and}}
	\qquad
		z \mapsto  \frac 1 {z - 2i},
	$$
and is depicted in Figure \ref{figr-apollonian}.

They key difference between $\cals$ and $\cala$ is that the former is  a hyperbolic attractor whereas the latter is parabolic: whereas the previous maps were contracting uniformly, here they are weakly attracting, having an indifferent fixed point. Visually, this corresponds to the fact that the largest circular holes approaching a vertex decay polynomially in size, not exponentially as in the previous example.

Correspondingly, the precise value for the Hausdorff and box-counting dimensions is unknown, but they do coincide. This is due to the conformality of the maps, as is the following formula of Boyd, which equates the so-called packing exponent to both Hausdorff and box-counting dimensions.
	\begin{prop}[Boyd \cite{boyd}]
		Enumerating the radii of all the circles excised from $\cala$ as $(r_k)_{k=1}^\infty$,
			$$
					\dim_B(\cala)
				=
					\dim_H(\cala)
				=
					\inf
						\bigg\{
							t > 0
						:
								\sum_{k=1}^\infty
									r_k^t
							<
								\infty
						\bigg\}.
			$$
	\end{prop}
The dimension has other useful formulations, for example thermodynamically in terms of the root of the topological pressure function \cite{fuchs,kantorovich-oh}. As for the value, Boyd \cite{boyd} showed $ \dim_H(\cala)\in [1.300197, 1.314534]$ and McMullen \cite{mcmullen} estimated $\dim_H(\cala) \approx 1.30568$, agreeing with more recent empirical estimates of Bai--Finch \cite{bai-finch} and others.

\subsubsection{The Rauzy gasket}

\begin{figure}[hbt]
	\centering
	\includegraphics[width=0.6\linewidth]{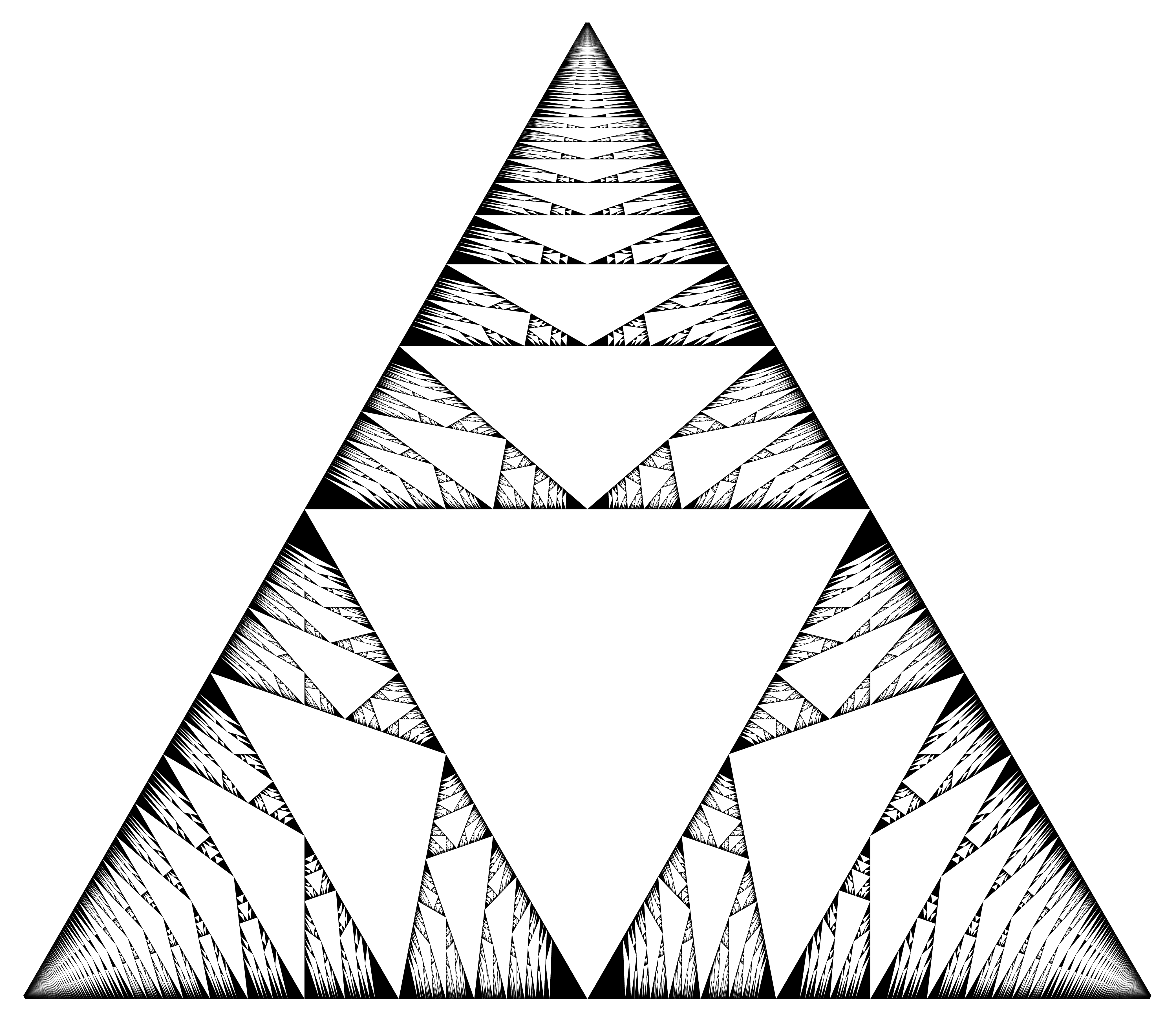}
	\caption[The Rauzy gasket $\calr$ (symmetric)]{The Rauzy gasket, $\calr$.}
	\label{fig:Rauzy gasket}
\end{figure}

The Rauzy gasket, $\calr$, in the fractal geometric sense represents the worst of all worlds. Not only is it a parabolic attractor like the Apollonian gasket, but the maps involved are non-conformal and highly distorting. 

Here, the three attracting maps $T_1,T_2,T_3$, acting on the standard 2-simplex, are the projective action of the matrices
	$$
			N_1
		=
			\left(
			\begin{array}{ccc}
				 1 & 1 & 1 \\
				 0 & 1 & 0 \\
				 0 & 0 & 1 \\
			\end{array}
			\right),
		\qquad
			N_2
		=
			\left(
			\begin{array}{ccc}
				 1 & 0 & 0 \\
				 1 & 1 & 1 \\
				 0 & 0 & 1 \\
			\end{array}
			\right),
		\qquad
			N_3
		=
			\left(
			\begin{array}{ccc}
				 1 & 0 & 0 \\
				 0 & 1 & 0 \\
				 1 & 1 & 1 \\
			\end{array}
			\right)
	$$
(general definitions are given below), and are explicitly given by
{
	\allowdisplaybreaks
	\begin{gather*}
		 	T_1(x,y,z)
		=
		 	\left(
		  		\frac{1}{2-x},
		  		\frac{y}{2-x},
		  		\frac{z}{2-x}
		  	\right),
	\\
		  	T_2(x,y,z)
		=
			  \left(
				  \frac{x}{2-y},
				  \frac{1}{2-y},
				  \frac{z}{2-y}
			  \right),
	\\
			T_3(x,y,z)
		=
	  		\left(
	  			\frac{x}{2-z},
	  			\frac{y}{2-z},
	  			\frac{1}{2-z}
	  		\right).
	\end{gather*}
}

Visually, the triangles in Figure \ref{fig:Rauzy gasket} decay polynomially in size as they approach the corners, and are often quite elongated. Correspondingly, much less is known in comparison to the other two examples, as we now describe.

\subsection[Previously known results]{Previous results about the area and dimension of $\calr$}

The Rauzy gasket has an interesting history, appearing for the first time in 1991 in the work of Arnoux and Rauzy \cite{arnoux-rauzy} in the context of interval exchange transformations, where it was conjectured that $\area(\calr)=0$. The gasket was rediscovered by Levitt in 1993, in a paper \cite{levitt} which also included a proof that $\area(\calr)=0$, due to Yoccoz. The gasket emerged for a third time in the work of De Leo and Dynnikov \cite{de leo-dynnikov}, this time in the context of Novikov's theory of magnetic induction on monocrystals.\footnote{See \cite{dynnikov-hubert-skripchenko} for the dichotomy between \cite{de leo-dynnikov} and \cite{arnoux-rauzy}.} They gave an alternative proof that $\area(\calr)=0$ and proposed the stronger result $\dim_H(\calr)<2$. Furthermore, a conjecture of Maltsev and Novikov \cite{maltsev-novikov} implies $ \dim_H(\calr) \in(1,2)$. But what is actually known?
	\begin{itemize}

	\item De Leo and Dynnikov in \cite{de leo-dynnikov} empirically estimate $\dim_B(\calr)\approx 1.72$.

	\item Avila, Hubert and Skripchenko in \cite{avila-hubert-skripchenko} prove $\dim_H(\calr) < 2$.

	\item Gut\'errio-Romo and  Matheus in \cite{gutierrez romo-matheus} show $\dim_H(\calr) \geq 1.19$.

	\item Fougeron in \cite{fou} shows $\dim_H(\calr) \leq 1.825$, using estimates of Baragar \cite{baragar}.

	\item Pollicott and the author in \cite{pollicott-sewell} prove $\dim_H(\calr) \leq 1.7415$.

\end{itemize}
Most recently, Jiao, Li, Pan and Xu in \cite{chinese} proved a formula expressing $\dim_H(\calr)$ as the critical exponent of a certain Poincar\'e series, as we will present and use explicitly later.

\subsection{Main results}

\subsubsection{Set-up}

We first set-up the context and notation needed to state our general results.

\begin{defn}
	Given an $(n+1)\times(n+1)$ matrix $N \in \RR^{n+1,n+1}$ such that all its entries are non-negative and $\det(N) = \pm1$, its \textbf{projectivisation} is the normalised map
		$$
				T:\De^{(n)} \to \De^{(n)},
			\qquad
				T(
					x
				)
				=
					\frac
						{N \cdot x}
						{|N\cdot x|}
				,
		$$
	where $|\cdot|$ denotes the $l^1$ norm and
		$
				\De^{(n)}
			=
				\big\{
					x \in [0,1]^{n+1}
				:
					|x| = 1
				\big\}
		$
	the standard $n$-simplex. (Note that $T$ is injective because $N$ is invertible.)
\end{defn}


Fix a dimension $d\geq 2$, and a collection of $(d+1)\times (d+1)$ matrices $N_1,\ldots,N_m$ as in the above definition. Denoting their respective projectivisations by $T_1,\ldots,T_m$ and $\De$ as $\De^{(d)}$, we take $\calg$ to be the limit set of these maps, i.e., the unique non-empty subset of $\De$ such that
	$$
			\calg
		=
			\bigcup_{j=1}^m
				T_j(\calg)
		.
	$$
This $\calg$ we call a \textbf{self-projective attractor}.

\begin{defn}
	We refer to a connected component of $\De \setminus \calg$ as a \textbf{hole} (of $\calg$), and a connected component of $\De\setminus\bigcup_{j=1}^m T_j(\De)$ as a \textbf{main hole}.
\end{defn}

Finally note that, when we refer to the interior or boundary of a set, we mean as a subset of the $d$-dimensional hyperplane $\{x\in\RR^{d+1}: |x| = 1\}\supset \De$.

\subsubsection{General results}

Our main theorem is the following.
\begin{thm}
	\label{thm:big theorem}
	Suppose that the images $T_1(\De),\ldots,T_m(\De)$ have disjoint interiors, that $\calg$ has zero $d$-dimensional Lebesgue measure, and that the main holes of $\calg$ are the interiors of $d$-simplices.
%
%
%
	Then, enumerating the holes as $\{\nab_n\}_{n=1}^{\infty}$,
%
%
	%
		\begin{align}
				\overline\dim_B(\calg)
			&=
				d
			+
				\inf
					\left\{
						t \geq -1
					:
						\sum_{n = 1}^\infty
							\vol(\nab_n)
							\In(\nab_n)^{t}
					<
						\infty
					\right\}
			,
			\label{eq:innies and outies}
		\end{align}
	where $\In(\nab_n)$ denotes the \textbf{inradius} of $\nab_n$, the radius of the largest ball contained in $\nab_n$, and $\vol(\nab_n)$ its $d$-dimensional volume.
\end{thm}
%
%

%
\begin{example}
	In the case of the Rauzy gasket, the main hole is the interior of the triangle with vertices $\big(\frac12,\frac12,0\big)$, $\big(0,\frac12,\frac12\big)$ and $\big(\frac12,0,\frac12\big)$, the $T_j(\De)$ meet each other only at their vertices, and $\area(\calr)=0$, as cited above. Thus Theorem \ref{thm:big theorem} applies.
\end{example}
\begin{figure}[htb]
\centering
\includegraphics[width = 0.9\linewidth]{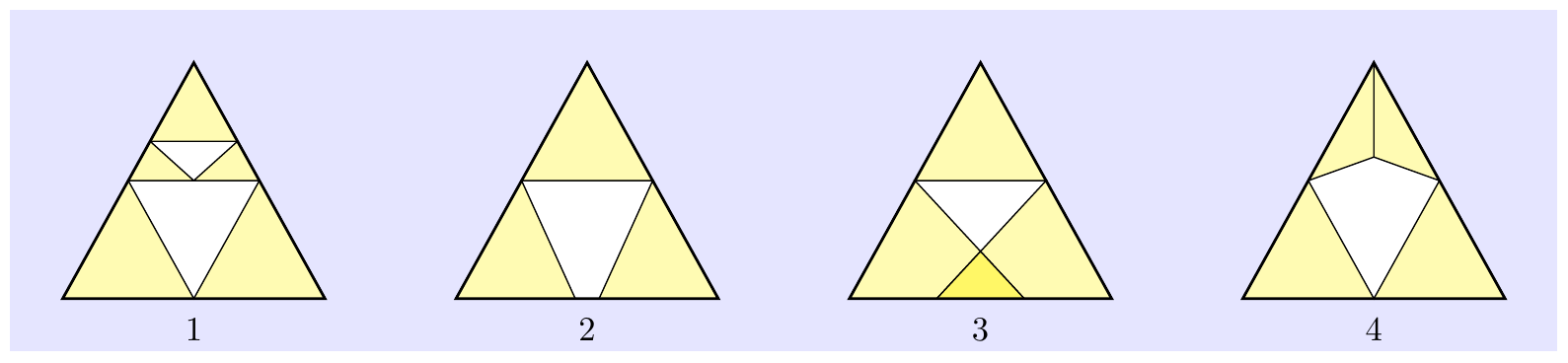}
\caption[little caption]
{Illustrating the geometric hypotheses of Theorem \ref{thm:big theorem}.}
\label{fig-badtriangles}
\end{figure}
\begin{example}
	To illustrate the geometric hypotheses of the theorem, Figure \ref{fig-badtriangles} shows four examples for $d=2$ (the images $T_j(\De)$ are shown in yellow and the main holes in white). Considering the hypotheses of Theorem \ref{thm:big theorem}, only in the first does the theorem apply:
	\begin{enumerate}
		\item The geometric hypotheses are satisfied here.
		\item The main hole contains part of the boundary of $\De$, so is not the interior of a triangle.
		\item The images do not have disjoint interiors.
		\item The main hole is the interior of a convex quadrilateral, not a triangle.
	\end{enumerate}
	Note that, in the last case, Theorem \ref{thm:convex theorem} below would give an upper bound on the dimension.
\end{example}

As mentioned in the example, the method of proof generalises to give an upper bound when the holes are merely convex, as follows.
\begin{thm}
	\label{thm:convex theorem}
	Suppose that $T_1(\De),\ldots,T_m(\De)$ have disjoint interiors, that
	$\vol(\calg)=0$,
	and that the main holes are the interiors of convex $d$-dimensional polytopes. Then, enumerating the holes by $\{\nab_n\}_{n=1}^{\infty}$,
		\begin{align}
				\overline\dim_B(\calg)
			&\leq
					d
			+
				\inf
					\left\{
						t \geq -1
					:
						\sum_{n = 1}^\infty
							\vol(\nab_n)
							\In(\nab_n)^{t}
					<
						\infty
					\right\}
			.
			\label{eq:convex}
		\end{align}
\end{thm}
Returning to Theorem \ref{thm:big theorem}, one of the main corollaries is the following result which partially confirms a conjecture of De Leo \cite{de-leo conjecture original}. It will be written in a more convenient way later, once some notation is established.
\begin{thm}
	Under the assumptions of Theorem \ref{thm:big theorem},
		\begin{align}
				\overline\dim_B(\calg)
			&\geq
				\frac{d-1}d
				\inf
				\left\{
					t > 0
				:
					\sum_{n=1}^\infty
					\sum_{|i|=n}
						\big\|
							N_{i_1}
							N_{i_2}
							\cdots 
							N_{i_n}
						\big\|^{-t}
				<
					\infty
				\right\}
			\\
			&=
				\frac{d-1}d
				\lim_{t \to \infty}
				\left(
					\frac1t
					\log
						\sum_{n=0}^\infty
					\#
						\bigg\{
							|i|=n
						:
							\|
								N_{i_1}
								N_{i_2}
								\cdots 
								N_{i_n}
							\|	
						<
							e^t
						\bigg\}
				\right)
			,
		\label{eq:lyapunov lower bound}
		\end{align}
	where $|i|=n$ is shorthand for $i \in \{1,2,\ldots,m\}^n$ and $\|\cdot\|$ is any matrix norm.
	\label{thm-de leo lower}
\end{thm}

\subsubsection{Application to the Rauzy gasket}

In the context of the Rauzy gasket, we combine Theorem \ref{thm:big theorem} with the main result of \cite{chinese} to show that the box-counting dimension of the Rauzy gasket exists and equals its Hausdorff dimension:

\begin{thm}
	$\dim_B(\calr) = \dim_H(\calr)$.
	\label{thm:equality}
\end{thm}

Regarding its value, along with the upper bound
		$
				\dim_H(\calr)
			\leq
				1.7415
		$%
, we give an effective lower bound via Theorem \ref{thm-de leo lower}, using estimates of \cite{sewell thesis} and the Mathematica code available at \cite{sewell github}:
\begin{thm}
		$$
				\dim_H(\calr)
			\geq
				1.6196
			.
		$$
	\label{thm:lower num}
\end{thm}

\subsection{Contents of this note}

In section 2, we introduce some necessary notation and definitions. In section 3, we relate the $\eps$-neighbourhood of $\calg$ to the interior neighbourhoods of the holes, and in section 4, we use the fact that the holes are simplicial to derive an explicit formula. In section 5 we deduce an expression for the Mellin--Laplace transform, which we use to complete the proof of Theorem \ref{thm:big theorem} in section 6. Section 7 is dedicated to proving the lower bound in Theorem \ref{thm-de leo lower}, and section 8 the upper bound in the case of convex holes. The final section 9 contains the applications to the Rauzy gasket, and the appendix contains the proofs of all our geometric lemmas.

\section{Definitions}

\subsection{Notation}

From now on, unless otherwise stated we fix the above notation and the assumptions of Theorem \ref{thm:big theorem}, to which we add the following convenient notation.

\begin{defn}
	Let
		$$
				\cali
			:=
				\{\emp\}
			\cup
				\bigcup_{n \in \NN}
					\{1,2,\ldots, m\}^n
		$$
 denote the set of finite words in $\{1,2,\ldots, m\}$. Given $i = (i_1,\ldots,i_n)\in\cali$, write
	\begin{itemize}

		\item $N_i = N_{i_1}\cdot N_{i_2}\cdot \ldots \cdot N_{i_n}$,

		\item $T_i = T_{i_1}\circ T_{i_2}\circ\cdots\circ T_{i_n}$, and

		\item $\De_i = T_i(\De)$.
	\end{itemize}
	Conventionally, $N_\emp$ refers to the identity matrix, $T_\emp$ to the identity map, and $\De_\emp$ to $\De$.
\end{defn}
Furthermore, since by assumption the main holes $\nab_1,\ldots,\nab_K$ are the interiors of $d$-simplices, one can construct matrices $M_1$,\dots,$M_K$ such that the closure of each $\nab_k$ is the image of the projectivisation, $S_k$ say, of $M_k$ (in the case of the Rauzy gasket, we omit the subscripts). Fixing these, we naturally extend our alphabet $\cali$ to $\cali^\ast = \cali\times\{1,2,\ldots,K\}$ and our notation as follows: if $i' = (i,k)\in \cali^\ast$, write
\begin{itemize}
	\item $N_{i'} :=N_iM_k$,
	\item $T_{i'}:=T_iS_k$, and
	\item $\nab_{i'} = T_i(\nab_k)$.
\end{itemize}
Note that our use of $\nab$ is to distinguish holes from images.

\subsection{Minkowski box-counting dimension}

	It is convenient here to use the following version of the box-counting dimension. All definitions and statements shall be made in the context of $\De = \De^{(d)}$ for simplicity.
\begin{defn}
	Given $\Om \subset \De$ closed and $\eps>0$, let
		$$
				\Om_\eps
			=
			\big\{
				x \in \De
			:
				\|x - y\|
			\leq
				\eps
			\text{ for some }
				y \in \Om
			\big\}
		$$
	denote the $\eps$-neighbourhood of $\Om$, where $\|\cdot\|$ is the standard Euclidean metric.
\end{defn}
\begin{defn}
	If $\calg \subset \De$ is closed, we define its upper box-counting dimension $\overline\dim_B(\calg)$ by
		$$
				\overline \dim_B
				(\calg)
			=
				d
			-
				\liminf_{\eps\to 0^+}
				\left(
					\frac
						{\log\big(\vol(\calg_\eps)\big)}
						{\log(\eps)}
				\right),
		$$
	and its lower box-counting dimension  $\underline\dim_B(\calg)$ by
		$$
				\underline \dim_B
				(\calg)
			=
				d
			-
				\limsup_{\eps\to 0^+}
				\left(
					\frac
						{\log\big(\vol(\calg_\eps)\big)}
						{\log(\eps)}
				\right).
		$$
	If the two limits coincide, its box-counting dimension $\dim_B(\calg)$ exists and takes this value.
\end{defn}

\section{Excised holes}

In this short section, we relate the $\eps$-neighbourhood of $\calg$ to the so-called $\eps$-inner neighbourhoods of the holes. This we do with two simple lemmas whose proofs can be safely omitted.

The first lemma formalises the intuitive, geometric description of $\calg$ as the limiting set of a process of excising holes, as fractals are often described informally.
More precisely, it asserts that any hole is the image of one of the main holes under the $T_i$, which follows simply from the standard definition of $\calg$ as a limit set:
	$$
			\calg
		=
			\bigcap_{n=1}^\infty
			\bigcup_{|i|=n}
				\De_i,
	$$
and also that $T_1,\ldots,T_m$ are injective and their images $\De_1,\ldots,\De_m$ have disjoint interiors.
\begin{lem}
	\label{lem:excised-holes}
	Every hole of $\calg$ is of the form $\nab_i$ for some $i \in \cali^\ast$. I.e.,
		\begin{equation}
				\De
			=
				\calg
			\cup
				\bigcup_{i \in \cali^*}
					\Int(\nab_i)
			,
		\label{eq:attractor as holes}
		\end{equation}
	where the unions are disjoint.
\end{lem}
%

The second lemma requires the following simple definition.
\begin{defn}
	For any $\Om \subset \De$ and $\eps>0$, define its \textbf{inner $\eps$-neighbourhood} $L_\eps(\Om)$ by
		 $$
		 		L_\eps(\Om)
		 	=
		 		\big\{
		 			x \in \Om
		 		:
		 			\|x - y\|
		 		\leq
		 			\eps
		 		\text{ for some }
		 			y \in \partial\Om
		 		\big\}
		 	,
		 $$
    the set of points in $\Om$ at most $\eps$ away from the boundary $\partial \Om$.
\end{defn}
This second lemma is a natural consequence of the first.
\begin{lem}
	For all $\eps>0$,
		\begin{equation}
				\calg_\eps
			=
				\calg
				\cup
				\bigcup_{i\in \cali^\ast}
					L_\eps(\nab_i).
		\end{equation}
		\label{lem:attractor-nhd-in-excised-form}
\end{lem}
%


\section{Simplices, inner neighbourhoods and inradii}

\begin{figure}[tbh]
	\centering
\begin{tikzpicture}

			\def\a{(0,0)}
			\def\b{(7,0)}
			\def\c{(6,3)}

		\fill[white,use as bounding box] (0,0) rectangle (8,4);

				\b--\c;

			\draw[line width=53,yellow!60,join=round] \a--\b--\c--cycle;
			\draw[line width=50,yellow!30,join=round] \a--\b--\c--cycle;

			\fill[white]
				(-2,0) -- (9,0) -- (9,-2) -- (-2,-2);
			\fill[white]
				(-2,-1) -- (8,4) -- (-2,4);
			\fill[white]
				\b -- \c -- (8,4) -- (8,0);

			\draw[thick] \a--\b--\c--cycle;

			\draw[->,shorten >= 1pt]  (5,0) -- (5,0.9) node[midway,anchor = west]{$\eps$};

			\draw (5,0) rectangle ++(0.15,0.15);

			\node at (3,2) {$\Om$};
			\node at (3,1/2) {$L_\eps(\Om)$};
\end{tikzpicture}

	\caption{The inner $\eps$-neighbourhood $L_\eps(\Om)$ of triangle $\Om$, when $\eps<\In(\Om)$.}
	\label{fig:inner-neighbourhood}
\end{figure}

In this section we make the inner neighbourhoods $L_\eps(\nab_i)$ explicit. This exploits the fact that an (injective) projective map on $\De$ maps any $k$-simplex onto another $k$-simplex. Hence, our assumption that the main holes are the interiors of $d$-simplices extends immediately to all holes, and we can apply the following lemmas. Their proofs, as with all geometric results, can be found in the appendix. In any case, we need a notion of inradius and incentre

The properties of the incentre of a simplex, described in the next lemma, are essential in the proofs of the two lemmas that follow.
\begin{lem}
	For any simplex $\Om$, there is a unique point $x \in \Om$ furthest from  $\partial\Om$, which we call the \textbf{incentre} of $\Om$.
\end{lem}
Our next result explicitly describes the inner neighbourhood of a simplex in terms of its incentre and inradius.
\begin{defn}
	For any convex body $\Om\subset \De$, let $\In(\Om)$ denote its \textbf{inradius}
		\begin{equation}
				\In(\Om)
			=
				\max\big\{
					\|x - y\|
				:
					x \in \Om,\
					y \in \partial\Om
				\big\}
		,
		\label{eq:inradius}
		\end{equation}
	i.e., the radius of the largest ball contained within $\Om$.
\end{defn}

\begin{lem}
	For any simplex $\Om$, $L_\eps(\Om)$ either equals $\Om$ or is obtained by excising a similar copy of it:
		\begin{equation}
				L_\eps(\Om)
			=
				\begin{cases}
					\Om,
				&
					\text{if }\eps \geq \In(\Om),
				\\
					\Om
				\setminus
					\Int\left[
						\left(
							1-
							 \frac
							 	{\eps}
							 	{\In(\Om)}
						\right)
					\cdot
						\Om
					\right],
				&
					\text{if }\eps < \In(\Om),
				\end{cases}
		\label{eq:L_epsilon-made-explicit}
		\end{equation}
	where, for $\lambda > 0$, $\la\cdot\Om$ denotes the similar copy to $\Om$ obtained by enlarging $\Om$ by a factor of $\lambda$ about its incentre.
	\label{lem:L_epsilon-made-explicit}
\end{lem}

\begin{rem}
	The second case of \eqref{eq:L_epsilon-made-explicit} is illustrated in Figure \ref{fig:inner-neighbourhood} when $\Om$ is a triangle.
\end{rem}

Our final result relates the inradius of a simplex to its volume and surface area as follows, generalising Heron's formula for the inradius of a triangle.

\begin{lem}
	For any non-degenerate $d$-simplex $\Om\subset \De$,
		$$
				\In(\Om)
			=
				d
			\cdot
				\frac
					{\vol(\Om)}
					{\per(\Om)}.
		$$
	where $\per(\Om)$ denotes the $(d-1)$-dimensional Lebesgue measure of the boundary $\partial\Om$.
	\label{lem:heron}
\end{lem}

\section[Laplace Transform]{Laplace transform of $\vol(\calg_\eps)$}


%
%

Our analysis begins after the following simple consequence of Lemmas  \ref{lem:attractor-nhd-in-excised-form} and \ref{lem:L_epsilon-made-explicit}.
%

%
\begin{prop}
	For any $d$-simplex $\Om \subset \De$,
		$$
				\vol(L_\eps(\Om))
			=
				\vol(\Om)
				\left(
					1
				-
					\max
					\left(
						0,
						1
					-
						\frac
							{\eps}
							{\In(\Om)}
					\right)^d
				\right)
	%
			.
		$$
	Consequently, since $\vol(\calg)=0$,
		\begin{align*}
					\vol(\calg_\eps)
			&=
				\sum_{i\in \cali^*}
					\vol\big(L_\eps(\nab_i)\big)
			=
				\sum_{i\in \cali^*}
					\vol(\nab_i)
					\left(
						1
					-
						\max
						\left(
							0,
							1
						-
							\frac \eps {\In(\nab_i)}
						\right)^d
					\right)
			.
		\end{align*}
	\label{prop:area-of-L_epsilon}
\end{prop}

\begin{figure}
	\centering
	\begin{tikzpicture}[scale=3]

				\fill[yellow!50] (3.5,1) -- (1,1) parabola (0,0) -- (3.5,0);

				\draw[dashed, gray] (0,0) rectangle (1,1);

				\draw[thick] (3.5,1) node[anchor = south east] {}
				-- (1,1) parabola (0,0);

				\draw[<->] (0,1.7) -- (0,0) -- (3.7,0); 

				\node at (1,-0.05) [anchor = north] {$\In(\nab_i)$};

				\node at (0,1) [anchor = east] {$\mathllap{\vol(\nab_i)}$};

	\end{tikzpicture}
	\caption{Sketch of $\eps \mapsto \displaystyle \vol(\nab_i)\left(1 - \max\left(0,1 - \frac \eps {\In(\nab_i)}\right)^d\right)$}
\end{figure}
The following lemma gives an explicit form for the Laplace transform
	$$
			f(t)
		:=
			\int_0^\infty
				e^{-tx}
				\vol(\calg_{\exp(-x)})
			\;\mathrm dx
		=
			\int_0^1
				\eps^{t-1}
				\vol(\calg_{\eps})
			\;\mathrm d\eps.
	$$
\begin{lem}
	For all $t > 0$,
		$$
				f(t)
			=
				\frac
					{\vol(\De)}
					t
			-
				\frac
					{d!}
					{t(t+1)(t+2)\cdots(t+d)}
				\sum_{i\in \cali^*}
					\vol(\nabla_i)
					\In(\nabla_i)^t
				.
		$$
	In particular, $f$ is analytic on the half-plane
		$
			\left\{
				z \in \CC
			:
				\Re(z)
			>
				\sig
			\right\}
		$%
	, where
		$$
				\sig
			:=
				\inf
					\left\{
							t \geq -1
						:
							\sum_{i\in \cali^*}
								\vol(\nab_i)
								\In(\nab_i)^t
						<
							\infty
					\right\}
			.
		$$
\end{lem}
\begin{proof}
	We prove the formula by direct computation, which we justify line-by-line afterwards.
{\allowdisplaybreaks
	\begin{align}
			f(t)
		&=
			\int_0^1
				\eps^{t-1}
				\vol(\calg_\eps)
			\;
			\mathrm d\eps
		\\
		&=
			\sum_{i\in \cali^\ast}
				\vol(\nab_i)
				\int_0^1
					\eps^{t-1}
					\left(
						1
					-
						\max
						\left(
							0,
							1
						-
							\frac \eps {\In(\nab_i)}
						\right)^d
					\right)
				\mathrm d\eps
		\\
		&=
			\sum_{i\in \cali^\ast}
				\vol(\nab_i)
				\left(
					\int_0^1
						\eps^{t-1}\;
					\mathrm d\eps
				-
					\int_0^{\In(\nab_i)}
						\eps^{t-1}
						\left(
							1
						-
							\frac \eps {\In(\nab_i)}
						\right)^d
					\mathrm d\eps
			\right)
		\\
		&=
			\sum_{i\in \cali^\ast}
				\vol(\nab_i)
				\left(
					\frac1t
				-
					\In(\nabla_i)^t
					\int_0^1
						y^{t-1}
						\left(
							1
						-
							y
						\right)^d
					\mathrm dy
			\right)
		\\
		&=
			\sum_{i\in \cali^\ast}
				\vol(\nab_i)
				\left(
					\frac1t
				-
					\frac
						{d!}
						{t(t+1)(t+2)\cdots(t+d)}
					\In(\nab_i)^t
			\right)
		\\
		&=
			\frac
				{\vol(\De)}
				t
		-
			\frac
				{d!}
				{t(t+1)(t+2)\cdots(t+d)}
			\sum_{i\in\cali^*}
				\vol(\nabla_i)
				\In(\nabla_i)^t
	,
	\end{align}
}
 as required. As for the explanation:
\begin{itemize}
	\item The first equality is the definition of $f$ given above.

	\item The second uses Proposition \ref{prop:area-of-L_epsilon} and exchanges integral and sum, using positivity.

	\item The third separates the two terms of the integral.

	\item The fourth uses the substitution $y = \eps/\In(\nabla_i)$.

	\item The fifth evaluates the Bernoulli integral
		$$
					\int_0^1
						y^{t-1}
						\left(
							1
						-
							y
						\right)^d
					\mathrm dy
				=
					\frac
						{d!}
						{t(t+1)(t+2)\cdots(t+d)}.
		$$
	\item The final equality simply expands the bracket, noting that
		$$
				\sum_{i\in \cali^*}
					\vol(\nabla_i)
			=
				\vol(\De),
		$$
	which is equivalent to $\vol(\calg) = 0.$
\end{itemize}

As for the analytic continuation, we note that, decreasing along the real line, the singularity we first meet can only occur at $z = 0,-1,\ldots,-d$, or at $\sig$. It therefore suffices to show that there is no pole at 0 itself, assuming $\sig\neq 0$: indeed, it is easy to see that the residues of the two terms at zero cancel perfectly in this case.
\end{proof}
%



\section{Proof of main theorem}

We are now ready to complete the prove of our main theorem, Theorem \ref{thm:big theorem}, which we restate in our new notation below.

\begin{thm}
	Under the assumptions of Theorem \ref{thm:big theorem}, $\overline\dim_B(\calg)=d+\sig$, i.e.,
		$$
				\overline\dim_B(\calg)
			=
				d
			+
				\inf\left\{
					t \geq -1
				:
					\sum_{i\in \cali^*}
	 					\vol(\nab_i)
	 					\In(\nab_i)^t
					<
						\infty\right\}
				.
		$$
\end{thm}

\begin{proof} We prove it in two inequalities: \\[5pt]
	\tikz[baseline=(A.base)]{
		\draw[use as bounding box] (0,0) node {$\geq$};
		\node[draw,fill=yellow!30,circle,inner sep=1pt] (A) {$\geq$};
	}
	Writing $\overline \dim_B(\calg) = r$, by definition,
		$$
			\liminf_{\eps\to 0^+}
				\frac
					{\log(\vol(\calg_\eps))}
					{\log(\eps)}
			=
				d-r.
		$$
	I.e., for any $\de > 0$, there exists $\eps'>0$ such that whenever $\eps<\eps'$,
		$$
				\frac
					{\log(\vol(\calg_\eps))}
					{\log(\eps)}
			\geq
				d - r - \de
		\qquad
		\iff
		\qquad
				\vol(\calg_\eps)
			\leq
				\eps^{d - r - \de}.
		$$
	Thus, whenever $t > r + \de - d$, trivial splitting and estimation shows $f(t)$ is finite:
		\begin{align*}
				f(t)
			&=
				\int_0^{\eps'}
					\eps^{t-1}
					\vol(\calg_{\eps})
				\;\mathrm d\eps
			+
				\int_{\eps'}^1
					\eps^{t-1}
					\vol(\calg_{\eps})
				\;\mathrm d\eps
			\\
			&\leq
					\int_0^{\eps'}
						\eps^{-1 + (t + d - r - \de)}
					\;
					\mathrm d \eps
				+
					\vol(\De)
					\int_{\eps'}^1
						\eps^{t-1}
					\;
					\mathrm d\eps
			<
				\infty.
		\end{align*}
	Thus $t \geq \sig$. The required inequality $\overline\dim_B(\calg)\geq d + \sig$ follows since $t$ and $\de$ are arbitrary.
	\\[3pt]
	\tikz[baseline=(A.base)]{
		\draw[use as bounding box] (0,0) node {$\leq$};
		\node[draw,fill=yellow!30,circle,inner sep=1pt] (A) {$\leq$};
	}
	We may assume $\sig < 0$, since otherwise $\calg \subset \De \implies \overline{\dim}_B(\calg) \leq d$ trivially suffices. Given any $\tau \in (\sig,0)$, the sum
		$$
				\sum_{i\in \cali^*}
 					\vol(\nab_i)
 					\In(\nab_i)^z
		$$
	is finite for $z = \tau$, and hence uniformly bounded on the line $\{\Re(z) = \tau\}$. Consider the Laplace transform
		$$
			\int_0^\infty
				\frac
					{\mathrm d}
					{\mathrm dx}
				\vol(\calg_{\exp(-x)})
				\cdot
				e^{-tx}
			\mathrm d x.
		$$
	Integrating by parts \cite[p.30]{davies}, this transform equals
		$$
				t\cdot\big(f(0^+) - f(t)\big)
			=
				\frac
					{d!}
					{(t+1)\cdots(t+d)}
					\sum_{i\in \cali^*}
						\vol(\nab_i)
						\In(\nab_i)^t
	 			.
		$$
	Since the right hand side is integrable on the line $\{\Re(z) = \tau\}$ and has no singularities on or to the right of it, the Laplace inversion theorem \cite[p.42ff]{davies} applies to give, for all $x>0$,
		\begin{align}
				\frac
					{\mathrm d}
					{\mathrm dx}
				\vol(\calg_{\exp(-x)})
			&=
				\frac 1 {2\pi i}
				\int_{\tau- i \infty}^{\tau + i \infty}
			\frac
				{d!}
				{(t+1)\cdots(t+d)}
				\sum_{i\in \cali^*}
					\vol(\nab_i)
					\In(\nab_i)^t
					\,
					e^{tx}
					\;
				\mathrm dt
			\\
			&=
				e^{\tau x}
				\int_{\tau- i \infty}^{\tau + i \infty}
					\calo\left(
						\frac
							1
							{(t+1)(t+2)\cdots(t+d)}
					\right)
				\mathrm dt
			\\
			&=
				\calo(e^{\tau x})
			.
		\end{align}
	In particular, since $\tau<0$, integrating this inequality gives
		$$
				\vol(\calg_{\exp(-x)})
			=
				\calo
					\left(
						\int_x^\infty
							e^{\tau y}
						\;\mathrm dy
					\right)
			=
				\calo\left(e^{\tau x}\right),
		$$
	i.e., $\vol(\calg_\eps) = \calo(\eps^{-\tau})$.	Thus, by definition, $\overline\dim_B(\calg) \leq d + \tau$, and we may take $\tau \to \sig$.
\end{proof}

\section{A lower bound for $\overline \dim_B(\calg)$}

In this section, we prove Theorem \ref{thm-de leo lower}, affirming the conjecture of De Leo \cite{de-leo conjecture original} under the geometric assumptions of our main theorem if one further assumes that $\dim_B\calg$ exists.

The new notation makes the definition of $\rho$ much easier to state, as we do below.
\begin{thm}
	Under the assumptions of Theorem \ref{thm:big theorem},
	\begin{align*}
			\overline \dim_B(\calg)
		\geq
			\max
			\left(
				d-1,
				\frac
					{d\cdot\rho}{d+1}
			\right),
	\end{align*}
	where
	\begin{align*}
			\rho
		&:=
			\inf\left\{
				r \geq 0
			:
					\sum_{i\in\cali}
						\|N_i\|^{-r}
				<
					\infty
			\right\}
		=
				\lim_{t\to \infty}
					\left(
						\frac
						{
							\log\#
							\{
								i \in \cali
							:
								\|N_i\|\leq t
							\}
						}
							{\log t}
					\right)
		.
	\end{align*}
\end{thm}

The proof of this theorem uses the following useful formula for the action of projective maps on the volume of simplices.

\begin{prop}
	Let $T:\De\to\De$ projectivise the matrix $N$ with $\det(N)=\pm1$. Then, for any $d$-simplex $\Om \subset \De$,
		$$
				\frac
					{\vol(T(\Om))}
					{\vol(\Om)}
			=
				\prod_{e:\text{ vertex of }\Om}
				|N \cdot e|^{-1}
				,
		$$
	where $|\cdot|$ denotes the $l^1$ norm.
	\label{prop-area form}
\end{prop}

\begin{proof}
	The case $\Om=\De$ was proven in \cite[Lem. 4.1]{pollicott-sewell-2}, and we show how this more general case follows: Let $M\in\text{SL}_d(\RR)$ be such that its projectivisation maps the vertices of $\De$ onto those of $\Om$. Then by linearity, the $\Om = \De$ case implies
		$$
				\frac
					{\vol(T(\Om_i))}
					{\vol(\Om)}
			=
				\prod_{e:\text{ vertex of }\De}
					\frac
						{|M\cdot e|}
						{|N\cdot M\cdot e|}
			=
				\prod_{e:\text{ vertex of }\De}
						\left|
							N\cdot
							\left(
								\frac
									{M\cdot e}
									{|M\cdot e|}
							\right)
						\right|^{-1}
			=
				\prod_{e:\text{ vertex of }\Om}
				|N \cdot e|^{-1},
		$$
	as required.
\end{proof}

\begin{proof}[Proof of Theorem \ref{thm:convex theorem}]
	We first prove the inequality involving the definition of $\rho$. For simplicity, we will use the symbol $\lesssim$ to denote that an inequality holds up to multiplication by a fixed constant, depending at most on $d$ and $s$.
	First observe that, for any $d$-simplex $\Om$, the volume of the largest contained ball cannot exceed the total volume:
		$
				\In(\nab_i)
			\lesssim
				\vol(\nab_i)^{1/d}.
		$
	Therefore, for any $t\in (-1,0)$,
		\begin{align*}
				\sum_{i\in \cali^*}
						\vol(\nab_i)
						\In(\nab_i)^t
			&\gtrsim
				\sum_{i\in \cali^*}
				\vol(\nab_i)^{1+t/d}
			\\
			&=
				\sum_{i \in \cali}
				\sum_{k=1}^K
					\vol(\nab_k)
				\prod_{e:\text{ vertex of }\nab_k}
					|N_i\cdot e|^{-(1+t/d)}
			\\
			&\geq
				\sum_{k=1}^K
					\vol(\nab_k)
				\cdot
				\sum_{i \in \cali}
					\|N_i\|^{-(d+1)(1+t/d)}.
		\end{align*}
	Thus, recalling $\sig$, the infimal value of $t\in(-1,0)$ for which the left hand side converges, it satisfies
		\begin{align*}
				0
			\geq
				\sig
			&\geq
				\inf\left\{
					t\geq -1
				:
					\sum_{i \in \cali}
						\|N_i\|^{-(d+1)(1+t/d)}
				<
					\infty
				\right\}
			=
				\max
				\left(
					-1
				,
					\frac
						{d\rho}
						{d+1}
				-
					d
				\right),
		\end{align*}
	and the inequality follows from Theorem \ref{thm:big theorem}.

	To show the equivalence of the two definitions of $\rho$, let $\rho$ denote the left hand side and $\la$ the right hand side, whose limit exists by submultiplicativity. For any $\eps>0$,
		\begin{align*}
				\sum_{i \in \cali}
					\|N_i\|^{-(\la+2\eps)}
			&=
				\sum_{n=1}^\infty
				\sum_{
					\substack{
						i \in \cali:\\
						2^{n-1}\leq\|N_i\|<2^{n}
					}
				}
					\|N_i\|^{-(\la+2\eps)}
			\\
			&\leq
				\sum_{n=1}^\infty
					2^{-n(\la+2\eps)}
					\#\{i \in \cali : \|N_i\|\leq 2^{-n}\}
			\\
			&\lesssim
				\sum_{n=1}^\infty
					2^{-n(\la+2\eps)}
					2^{n(\la+\eps)}
			<
				\infty.
		\end{align*}
	Hence $\la \geq \rho$. The reverse inequality is similar, and completes the proof.
	%
	%
%
\end{proof}

\section{An upper bound for $\overline\dim_B(\calg)$ in the convex hole case}

In this penultimate section we outline how, relaxing the assumptions so that the main holes may be convex polytopes, the proof of Theorem \ref{thm:big theorem} can be extended to give an upper bound. Previously given as Theorem \ref{thm:convex theorem}, it can be restated in our new notation as follows.
\begin{thm}
	Under the assumptions of Theorem \ref{thm:convex theorem},
		$$
				\overline\dim_B(\calg)
			\leq
				d
			+
				\inf\left\{
					t \geq -1
				:
					\sum_{i\in \cali^*}
	 					\vol(\nab_i)
	 					\In(\nab_i)^{t}
					<
						\infty\right\}
				.
		$$
\end{thm}
The fundamental observation in the proof is that projectivised linear maps preserve convexity (so that the image of a convex polytope is a convex polytope), which can be shown directly from the definition.
The main obstacle is that inner neighbourhoods of general convex bodies are more difficult to work with; we nevertheless obtain a general inclusion and inequality.

Note that, since there may be multiple points furthest from the boundary of a convex body $\Om$ (e.g., when $\Om$ the rectangle), we will denote any one of these points as an incentre.
\begin{lem}
	For any convex body $\Om\subset \De$,
		\begin{equation}
				L_\eps(\Om)
			\subset
				\begin{cases}
					\Om,
				&
					\text{if }\eps \geq \In(\Om),
				\\
					\Om
				\setminus
					\Int\left[
						\left(
							1-
							 \frac
							 	{\eps}
							 	{\In(\Om)}
						\right)
					\cdot
						\Om
					\right],
				&
					\text{if }\eps < \In(\Om);
				\end{cases}
		\label{eq:L_epsilon-made-explicit convex verzio}
		\end{equation}
	where, for $\lambda > 0$, $\la\cdot\Om$ denotes the similar copy to $\Om$ obtained by enlarging $\Om$ by a factor of $\lambda$ about any fixed incentre.
	Consequently,
		$$
				\vol(L_\eps(\Om))
			\leq
				\vol(\Om)
				\left[
					1
				-
					\max\left(
						0
					,
						1
					-
						\frac
							\eps
							{\In(\Om)}
					\right)^d
				\right]
		$$
	and
		\begin{equation}
				\vol(\calg_\eps)
			\leq
				\sum_{i\in\cali^\ast}
					\vol(\nab_i)
				\left[
					1
				-
					\max\left(
						0
					,
						1
					-
						\frac
							\eps
							{\In(\nab_i)}
					\right)^d
				\right]
			.
		\label{eq:inequooo}
		\end{equation}
\end{lem}


\begin{proof}
	We defer the proof of the first inclusion to the appendix, and simply note that the inequalities are immediate, the second following by combining the first with \ref{lem:attractor-nhd-in-excised-form} and the knowledge that all holes $\nab_i$ are convex and the assumption that $\vol(\calg)=0$.
\end{proof}

From this point, the proof of the theorem continues analogously to that of Theorem \ref{thm:big theorem}, applying the results of section 5 and onwards to the right hand side of \eqref{eq:inequooo}.

Finally, for completeness, we generalise the formula for the inradius of a simplex to an inequality. It is proven in the appendix.
\begin{prop}
	For any non-empty convex set $\Om\subset\De$,
		$$
				\frac
					{\vol(\Om)}
					{\per(\Om)}
			\leq
				\In(\Om)
			\leq
				(d-1)
				\frac
					{\vol(\Om)}
					{\per(\Om)}.
		$$
\end{prop}

\section{Application to the Rauzy gasket}

In this final section, we prove the important results on the equality of Hausdorff and box-counting dimensions of the Rauzy gasket and go on to outline the proof of the numerical lower bound provided by Theorem \ref{thm-de leo lower}

In this section, given a matrix $N$, we will denote by $\sig_j$ its $j$th largest singular value. The following result of \cite{chinese} is reminiscent of singularity dimension, first introduced by Falconer.
	\begin{thm}[\cite{chinese}]
		The Hausdorff dimension of the Rauzy gasket $\dim_H(\calr)$ is equal to the infimal value of $s\in(1,2)$ for which the Poincar\'e series
			\begin{align}
					\sum_{i \in \cali}
						\left(
							\frac
								{\sig_2}
								{\sig_1}
						\right)
						\left(
							\frac
								{\sig_3}
								{\sig_1}
						\right)^{s-2}
						(N_i)
			\label{eq:lambda}
			\end{align}
		converges.
	\end{thm}
Note that we will sometimes, for convenience, denote the summand of \eqref{eq:lambda} by $\phi_s(N_i)$.

We will apply this result by showing the following.
	\begin{lem} For any $s \in (1,2)$,
		$$
				\sum_{i \in \cali^\ast}
					\area(\nab_i)
					\In(\nab_i)^{s-2}
			\lesssim
				\sum_{i \in \cali}
					\left(
						\frac
							{\sig_2}
							{\sig_1}
					\right)
					\left(
						\frac
							{\sig_3}
							{\sig_1}
					\right)^{s-2}
					(N_i).
			$$
		%
		\label{lem:series}
	\end{lem}
Thus combining this lemma with Theorem \ref{thm:big theorem} and the above result, we see that $\overline\dim_B(\calr)\leq \dim_H(\calr)$, immediately giving Theorem \ref{thm:equality}.

\begin{proof}[Proof of Lemma \ref{lem:series}]

Following \cite{chinese}, we note that the Euclidean distance between $x,y \in \De$, $\|x-y\|$ equals, up to a constant factor, the area of the triangle formed by $x,y$ and the origin:
	$$
			\|x-y\|
		\asymp
			\|x\wedge y\|,
	$$
where $\asymp$ denotes that the two sides are comparable, i.e. their ratio is bounded away from 0 and $\infty$, by constants depending at most on $d$ and $s$.
Heron's formula (Lemma \ref{lem:heron}) gives, for $s \in(1,2)$,
	\begin{align}
			\sum_{i \in \cali^\ast}
				\area(\nab_i)
				\In(\nab_i)^{s-2}
		\asymp
			\sum_{i \in \cali^\ast}
				\area(\nab_i)
				\left(
					\frac
						{\per(\nab_i)}
						{\area(\nab_i)}
				\right)^{2-s}
			.
	\label{eq:seriesagain}
	\end{align}
Now, given $i \in\cali^*$, applying the area formula given in Proposition \ref{prop-area form},
	\begin{align*}
					\frac
						{\per(\nab_i)}
						{\area(\nab_i)}
				&\asymp
					\sum_{\{j,k,l\}= \{1,2,3\}}
						\|T_i(e_j)\wedge T_i(e_k)\|\cdot |N_i\cdot e_1||N_i\cdot e_2||N_i\cdot e_3|
				\\
				&=
					\sum_{\{j,k,l\}= \{1,2,3\}}
						\|N_i\cdot e_j\wedge N_i\cdot e_k\|\cdot |N_i\cdot  e_l|
				\\
				&\lesssim
					\phantom{\sum_{\{j,k,l\}= \{1,2,3\}}}
					3 \qquad \sig_1\sig_2\qquad\!\!\!\cdot\qquad\!\!\!\sig_1(N_i)
				\ =\ \sig_1^2\sig_2(N_i).
	\end{align*}
Hence, substituting into \eqref{eq:seriesagain},
	\begin{align}
			\sum_{i \in \cali^\ast}
				\area(\nab_i)
				\In(\nab_i)^{s-2}
		\lesssim
			\sum_{i \in \cali^\ast}
				\area(\nab_i)\cdot
				\sig_1^{4-2s}\sig_2^{2-s}(N_i)
			.
		\label{eq:stanley}
	\end{align}
To complete the proof, it will suffice to show $\area(\nab_i)\lesssim \sig_1^{-3}(N_i)$ for all $i\in \cali^\ast$, since then the right hand side of \eqref{eq:stanley} will be dominated by
	\begin{equation}
			\sum_{i \in \cali^\ast}
				\sig_1^{-3}
				\sig_1^{4-2s}\sig_2^{2-s}\,(N_i)
		=
			\sum_{i \in \cali^\ast}
				\left(
					\frac
						{\sig_2}
						{\sig_1}
				\right)
				\left(
					\frac
						{\sig_3}
						{\sig_1}
				\right)^{s-1}
				(N_i)
		\asymp
			\sum_{i' \in \cali}
				\left(
					\frac
						{\sig_2}
						{\sig_1}
				\right)
				\left(
					\frac
						{\sig_3}
						{\sig_1}
				\right)^{s-1}
				(N_{i'}),
		\label{eq:two sums compared}
	\end{equation}
where the first inequality uses that $\det(N_i) = 1 = \sig_1\sig_2\sig_3(N_i)$ for any $i\in\cali^\ast$ and the second relation is justified by term-by-term comparison: for any $i' \in \cali$, we have, for its corresponding $i \in \cali^\ast$,
	$
		N_i = N_{i'}\cdot M,
	$
where $M$ is a non-singular matrix to be made explicit below. Since $M$ is non-singular and $\det(N_i) = \det(N_{i'})=1$,
	$$
		\sig_1(N_i)\asymp \sig_1(N_{i'}),
		\qquad
		\sig_3(N_i)\asymp \sig_3(N_{i'})
		\implies
		\sig_2(N_i) = \frac 1 {\sig_1\sig_3}(N_i) \asymp  \frac 1 {\sig_1\sig_3}(N_i) = \sig_2(N_{i'})
	$$
and hence $\phi(N_i)\asymp \phi (N_{i'})$, showing that the two series in \eqref{eq:two sums compared} are comparable.

We prove the assertion  $\area(\nab_i)\lesssim \sig_1^{-3}(N_i)$ in three cases, writing $N = N_i$ for simplicity:\\[8pt]
\noindent \textbf{Case 1:} $i = (i',\cdot)\in\cali^*$, where $i'\in\cali$ contains all three symbols.

	Without loss of generality, write $N = A\cdot B$, where
		$$
				B
			=
				N_1\cdot N_2^a\cdot N_3^b\cdot M
		$$
	for some $a,b\in\NN$, fixing $M$ to be the matrix with the following explicit form:
		$$
				M
			=
				\begin{pmatrix}
					0&\alf&\alf\\
					\alf&0&\alf\\
					\alf&\alf&0
				\end{pmatrix},
		$$
	where $\alf = 2^{-1/3}$, so that $\det(M) = 1$. More explicitly,
		$$
				B
			=
				\alf
				\left(
					\begin{array}{ccc}
						 a b+a+b+2 & a b+2 a+b+2 & 2 a b+a+2 b+2 \\
						 a b+a+1 & a b+2 a & 2 a b+a+1 \\
						 b+1 & b+1 & 2 b \\
					\end{array}
				\right).
		$$
	We can see that the ratio of any two elements in a row of $B$ is at most two:
		$$
			\max_{j,k,l}
				\left(
					\frac
						{B_{j,k}}
						{B_{j,l}}
				\right)
			\leq
				2,
		$$
	independently of $a$ and $b$. Applying this directly gives, for	$\la = (1,1,1) \cdot A$ gives
		$$
				\frac
					{|N \cdot e_k|}
					{|N \cdot e_l|}
			=
				\frac
					{\la \cdot B \cdot e_l}
					{\la \cdot B \cdot e_k}
			=
				\frac
					{
					\sum_{j=1}^3
						\la_j B_{j,l}
					}
					{
					\sum_{j=1}^3
						\la_j B_{j,k}
					}
			\leq
				\frac
					{
					\sum_{j=1}^3
						\la_j (2B_{j,k})
					}
					{\phantom2
					\sum_{j=1}^3
						\la_j B_{j,k}
					}
			\leq
				2.
		$$
	Applying this in the area formula (Proposition \ref{prop-area form}) finishes the proof in this case:
		\begin{align}
				\area(\nab_i)
			\asymp
				\prod_{j=1}^3
					|N\cdot e_j|^{-1}
			\leq
				4
					\big(
					\max_j|N\cdot e_j|
					\big)^{-3}
			=
				4\cdot\|N\|_{l^1\to l^1}^{-3}
			\asymp
				\sig_1^{-3}(N),
		\label{eq:area domination end}
		\end{align}
	where $\|N_i\|_{l^1\to l^1}$ denotes the $l^1$ operator norm of $N_i$, using equivalence of norms at the end.\\[8pt]
\noindent \textbf{Case 2:} $i = (i',\cdot)\in\cali^*$; $i'\in \cali$ contains precisely two symbols.

	We assume symbols 1 and 2 appear in $i'$ (in that order), and claim that $N$ is of the form
			$$
					\alf
					\left(
						\begin{array}{ccc}
							 x_1 & x_2 & x_3 \\
							 y_1 & y_2 & y_3 \\
							 1 & 1 & 0 \\
						\end{array}
					\right),
			$$
		where $x_j, y_j$ are positive integers such that
			$$
				\max_{j,k}
				\left(
					\frac
						{x_j}
						{x_k}
				\right)<3
				\qquad\text{and}\qquad
				\max_{j,k}
				\left(
					\frac
						{y_j}
						{y_k}
				\right)
				<3.
			$$
		This follows by a simple induction. As a base case, for $b \in \NN$ consider the matrix
			$$
					N_1N_2^bM
				=
					\alf
					\left(
						\begin{array}{ccc}
							 b+2 & 2 b+1 & b+1 \\
							 b+1 & 2 b & b+1 \\
							 1 & 1 & 0 \\
						\end{array}
					\right),
			$$
		which is of the required form. Inductively assuming that $N$ is of the above form,
			\begin{align*}
					N_1N
				&=
					\alf\left(
						\begin{array}{ccc}
							 x_1 + y_1 + 1 & x_2 + y_2 + 1 & x_3 + y_3 \\
							 y_1 & y_2 & y_3 \\
							 1 & 1 & 0 \\
						\end{array}
					\right),
				\\
					N_2N
				&=
					\alf\left(
						\begin{array}{ccc}
							 x_1 & x_2 & x_3 \\
							 x_1 + y_1 + 1 & x_2 + y_2 + 1 & x_3 + y_3 \\
							 1 & 1 & 0 \\
						\end{array}
					\right),
			\end{align*}
		which again are easily shown to be of the required form, proving the claim.

		Thus we may write, for any column indices $l$ and $k$,
		$$
			\frac
				{|N \cdot e_l|}
				{|N \cdot e_k|}
		=
			\frac
				{
					x_l + y_l
				+
				 	\de_{\{l\neq 3\}}
				}
				{	x_k + y_k
				+
					\de_{\{k\neq3\}}
				}
		\leq
			\frac
				{
					(3 x_k -1 )+ (3y_k - 1)
				+
				 	1
				}
				{x_k + y_k
				}
		<
			3.
		$$
Hence an analogue of \eqref{eq:area domination end}  applies to show that $\area(\nab_i) \lesssim \sig_1^{-3}(N)$ (with a larger implied constant).\\[8pt]
\noindent \textbf{Case 3:} $i = (i',\cdot)\in\cali^*$; $i'\in \cali$ contains only one symbol.
	Assuming without loss that $i' = (1,\ldots,1)\in \{1\}^b$ for $b \in\NN_0$, we explicitly write
		$$
				N
			=
				N_1^bM
			=
				\alf
				\left(
					\begin{array}{ccc}
						 2 b & b+1 & b+1 \\
						 1 & 0 & 1 \\
						 1 & 1 & 0 \\
					\end{array}
				\right)
			.
		$$
	Hence $\displaystyle \max_{k,l}\left(\frac {|N \cdot e_k|}{|N \cdot e_l|}\right) = \frac {2b+2}{b+2} \leq 2$ and \eqref{eq:area domination end} applies as before to complete the proof.
\end{proof}

In order to outline the proof of Theorem \ref{thm:lower num}, note that the above proof, particularly \eqref{eq:area domination end}, can be altered to show more strongly that $\area(\De_{i}) \asymp \area(\nab_{i'}) \asymp \|N_i\|^{-3}$ in the above notation.

This observation allows us to relate the theoretical lower bound in Theorem \ref{thm-de leo lower} to a numerical lower bound as per the following result. In particular, the author in \cite{sewell thesis} uses renewal theory and Mathematica to give the upper and lower bounds. See \cite{sewell github} for the code used and an in-depth explanation thereof. This altogether gives the following.
\begin{thm}[\cite{sewell thesis}]
	$$
			\inf\left\{
				\de>0
			:
				\sum_{i \in \cali}
					\|N_i\|^{-3\de}
			<
				\infty
			\right\}
		=
			\inf\left\{
				\de>0
			:
				\sum_{i \in \cali}
					\area(\De_i)^{\de}
			<
				\infty
			\right\}
		\in
			[0.8098, 0.8204].
	$$
\end{thm}
Hence, by Theorems \ref{thm-de leo lower} and \ref{thm:equality}, $\dim_H(\calr) \geq 2 \times 0.8098 = 1.6196$, as claimed.


\appendix

\section{Geometric results}

In this appendix, we give elementary and complete proofs of our geometric results in two parts. The first concerns simplices, and the second some generalisations to complex bodies. We will state and prove them in the context of $\RR^n$ for simplicity.

\subsection{Inradii and incentres of simplices}

In this subsection, let $\Om\subset\RR^n$ denote a given non-degenerate $n$-simplex.

One key property of simplices (compared to general convex bodies) is that the incentre is uniquely defined, as per the following simple but important lemma.

\begin{lem}
	There is a unique point in $\Om$ furthest from the boundary $\partial\Om$ (the incentre), which is equidistant from all faces of $\Om$.
\end{lem}
\begin{figure}
	\centering
	\begin{tikzpicture}[scale=3]

		\path
			coordinate (O) at (0,0)
			coordinate (up) at (20:2)
			coordinate (down) at (-20:2);

		\foreach \r/\d in {0.345/1.395}
		{
			\draw[dashed,fill=blue!10](\d,0) circle (\d*\r);
			\draw[fill = yellow!80] (1,0) circle (\r) node {$B$};
			\draw[<->,shorten <= 1pt,shorten >= 1pt]
				let
					\p1 = (down)
				in
					(1+\r,0) -- (\x1,0)
				node[midway,anchor = south]
					{$d(B,S)$};
		}

		\draw[<->,shorten <= 1pt,shorten >= 1pt]
			let
				\p1 = (down)
			in
				(0,\y1-3pt) -- (\x1,\y1-3pt)
			node[midway,anchor = south]
				{$d(e,S)$};

		\draw[black!30,dotted]
			let
				\p1 = (down)
			in
				(0,0) -- (0,\y1-3pt)
				(\x1,0) -- (\x1,\y1-3pt)
			;

		\draw[thick]
			(O)
				node [anchor = east] {$e$}
		-- (up)
		-- (down)
			node [midway,anchor = west] {$S$}
		-- cycle;

	\end{tikzpicture}
	\caption{Enlarging ball $B \subset \Om$ to meet side $S$.}
	\label{fig:enlarge-your-balls}
\end{figure}
\begin{proof}
	Let $B$ be a ball contained in $\Om$, and suppose it does not meet side $S\subset\partial \Om$. Then, letting $e$ denote the opposing vertex, one can show the enlargement of $B$ about $e$ with scale
		$$
			\frac
				{d(e,S)}
				{d(e,S)-d(B,S)}:
		$$
	see Figure \ref{fig:enlarge-your-balls}, is larger and contained in $\Om$. Thus, a maximal ball in $\Om$ touches all of its sides, i.e., any point which is furthest from $\partial\Om$ is equidistant from all sides.

	We further claim there can be only one point equidistant from all sides: Since $\Om$ is non-degenerate, the normal vectors to each face are linearly independent: thus, if there were two such distinct points, moving from one to the other would increase distance to some side and decrease distance to another, a contradiction.
\end{proof}

The incentre provides a convenient decomposition of $\Om$ which we use in the next two proofs. We begin with the simpler of the two; the formula for the inradius.

\begin{lem}
		$$
				\In(\Om)
			=
				n
				\frac
					{\vol(\Om)}
					{\per(\Om)}.
		$$
\end{lem}

Thanks are due to Andr\'as Sandor, who produced the proof idea in about a minute!

\begin{proof}
	For a given side $S\subset \partial \Om$, let $\Om_S$ denote the set of points in $\Om$ whose closest side is $S$.
	We know that $\Om_S$ is bounded between $S$ itself and each of the hyperplanes bisecting $S$ and one of the other faces --- i.e., the set of points equidistant to these two sides. By the previous lemma, these hyperplanes must all meet at the incentre of $\Om$, thus $\Om_S$ is also an $n$-simplex. Taking $S$ as the base of $\Om_S$, it has height $\In(\Om)$, and hence volume
		$$
				\vol(\Om_S)
			=
				\frac 1n
				\vol_{n-1}(S) \cdot
				\In(\Om)
			.
		$$
	The formula follows by summing over $S$.
\end{proof}

Next we prove the important expression for the inner $\eps$-neighbourhood.

\begin{lem}
	For all $\eps>0$,
		$$
				L_\eps(\Om)
			=
				\begin{cases}
					\Om,
				&
					\text{if }
					\eps \geq \In(\Om),
				\\
					\Om
				\setminus
					\Int
					\left[
						\left(
							1
						-
							\frac
								\eps
								{\In(\Om)}
						\right)
					\cdot
						\Om
					\right],
				&
					\text{if }
					\eps \leq \In(\Om);
				\end{cases}
		$$
	where $x \mapsto \la\cdot x$ denotes enlargement, scale factor $\la$, about the incentre of $\Om$.
\end{lem}
\begin{figure}
	\centering
	\begin{tikzpicture}[scale=2]

		\path
			coordinate (O) at (0,0)
			coordinate (up) at (-115:2)
			coordinate (down) at (-65:2)
			coordinate (lup) at (-115:1.4)
			coordinate (ldown) at (-65:1.4)
			coordinate (julia) at (-75:1.6)
			coordinate (kisjulia) at (-75:{1.6*1.4/2})
			;

		\draw[dotted]
			let
				\p1 = (julia),
				\p2 = (up),
				\p3 = (down)
			in
				(\x2-27,\y1) -- (0,\y1)
			;
		\draw[dotted]
			let
				\p1 = (julia),
				\p2 = (up),
				\p3 = (down)
			in
				(\x2-27,0) -- (\x3+7,0)
			;

		\draw[dotted]
			let
				\p1 = (kisjulia),
				\p2 = (up),
				\p3 = (down)
			in
				(\x3+7,\y1) -- (0,\y1)
			;

		\draw[dotted]
			let
				\p1 = (up),
				\p2 = (down)
			in
				(\x1-27,\y1) -- (\x2+7,\y1)
			;

		\path[fill=white]
			(up) -- (down) -- (O);

		\draw[fill = yellow!50]
			(lup) -- (ldown)
			-- (O)
			;

		\node at (julia) [fill=blue!90, circle,inner sep = 1pt] {};

		\node at (julia) [anchor = south]{$\;x$};

		\node at (kisjulia) [anchor = south]{$\la\cdot x\quad$};

		\node at (kisjulia) [fill=blue!,circle,inner sep = 0.7pt] {};

		\path
			let
				\p1 = (down)
			in
		coordinate (ono) at (O);

		\node at (ono) [anchor = south] {$\phantom{{}_S}\Om_S$};


		\draw[shorten <= 3pt, shorten >= 3pt, white, line width =5pt]
			let
				\p1 = (up)
			in
				(\x1 - 7pt,\y1) -- (\x1-7pt,0);

		\draw[<->]
			let
				\p1 = (up)
			in
				(\x1 - 7pt,\y1) -- (\x1-7pt,0)
					node[midway,anchor = east]
					{$\In(\Om)$}
				;

		\draw[<->]
			let
				\p1 = (down),
				\p2 = (kisjulia)
			in
				(\x1 + 5pt,0) -- (\x1 + 5pt,\y2)
					node[midway,anchor = west]
					{$\In(\Om) - d_S(\la \cdot x)$};

		\draw[<->]
			let
				\p1 = (down),
				\p2 = (kisjulia)
			in
				(\x1 + 5pt, \y2) -- (\x1 + 5pt,\y1)
				node[midway,anchor = west]
					{$d_S(\la\cdot x)$}
				;

		\draw[<->]
			let
				\p1 = (up),
				\p2 = (julia)
			in
				(\x1 - 25pt, \y2) -- (\x1 - 25pt,\y1)
				node[midway,anchor = east]
					{$d_S(x)$}
				;

		\draw[<->]
			let
				\p1 = (up),
				\p2 = (julia)
			in
				(\x1 - 25pt, \y2) -- (\x1 - 25pt,0)
				node[midway,anchor = east]
					{$\In(\Om) - d_S(x)$}
				;


		\draw[thick]
			(O)  --
			(up) --
			(down)
				node[midway,anchor = north]
					{$S$}
			--cycle;

	\end{tikzpicture}

	\caption{The action of $x \mapsto \la\cdot x$ on $\Om_S$, for $\la \in (0,1)$.}

	\label{fig:lambda contraction}
\end{figure}
\begin{proof}
	It suffices to show, given $\eps < \In(\Om)$, that
		\begin{equation}
				\big\{
					x \in \Om
				\mid
						d(x,\partial\Om)
					\geq
						\eps
				\big\}
			=
				\left(
					1
				-
					\frac
						\eps
						{\In(\Om)}
				\right)
				\cdot \Om
		\label{eq:inner neighbourhood 1}
		\end{equation}
	and
		\begin{equation}
				\big\{
					x \in \Om
				\mid
						d(x,\partial\Om)
					=
						\eps
				\big\}
			=
				\left(
					1
				-
					\frac
						\eps
						{\In(\Om)}
				\right)
				\cdot \partial\Om.
		\label{eq:inner neighbourhood 2}
		\end{equation}
	Moreover, recalling the notation from the previous proof,
		$
				x
			\mapsto
				\left(
					1
				-
					\frac
						\eps
						{\In(\Om)}
				\right)
				\cdot
				x
		$
	maps each $\Om_S$ into itself, since it is a linear contraction towards one of the vertices of $\Om_S$. Hence, we only need to prove \eqref{eq:inner neighbourhood 1} and \eqref{eq:inner neighbourhood 2} with $\Om_S$ in place of $\Om$ and $S$ in place of $\partial \Om$: i.e.,
		$$
				\big\{
					x \in \Om_S
				\mid
						d_S(x)
					\geq
						\eps
				\big\}
			=
				\left(
					1
				-
					\frac
						\eps
						{\In(\Om)}
				\right)
				\cdot \Om_S,
		$$
	where $d_S := d(x,S)$ for $x \in \Om$. Indeed, by linearity, for any $x \in \Om_S$,
		$$
					d_S
					\left(
						\left(
							1
						-
							\frac
								\eps
								{\In(\Om)}
						\right)
					\cdot
						x
					\right)
				-
					\In(\Om)
			=
				\left(
					1
				-
					\frac
						\eps
						{\In(\Om)}
				\right)
				\big(
					d_S(x)
				-
					\In(\Om)
				\big)
		$$
	(see Figure \ref{fig:lambda contraction}), which rearranges to
		\begin{equation}
				d_S
					\left(
						\left(
							1
						-
							\frac
								\eps
								{\In(\Om)}
						\right)
					\cdot
						x
					\right)
				=
					\left(
						1
					-
						\frac
							\eps
							{\In(\Om)}
					\right)
					d_S(x)
				+
					\eps
				\geq
					\eps,
		\label{eq:contraction inequality}
		\end{equation}
	with equality iff $x \in S$.
	Moreover, we may extend $d_S$ linearly to the cone $C_S$ based at the incentre and containing $S$,\footnote{Equivalently, $C_S$ can be defined as those points $x \in \RR^n$ whose closest point in $\partial \Om$ lies on $S$.}
	so that	$d_S(x) < 0$ for $x \in C_S \setminus \Om_S$. The equality in \eqref{eq:contraction inequality} extends to $x \in C_S$ by linearity, and thus shows, for $x \in C_S$,
		$$
				d_S
					\left(
						\left(
							1
						-
							\frac
								\eps
								{\In(\Om)}
						\right)
					\cdot
						x
					\right)
			\geq
				\eps
			\qquad
			\text{if and only if}
			\qquad
				d_S(x)\geq 0
			\iff 
				x\in \Om.
		$$
	Taking the union over sides $S$ gives \eqref{eq:inner neighbourhood 1} and \eqref{eq:inner neighbourhood 2}; thus the lemma follows.
\end{proof}

\subsection{Inradii and incentres of convex bodies}

We now fix $\Om \in \RR^n$ to be a convex body (i.e., a closed convex set with non-empty interior).

We first give the proof of the inclusion \eqref{eq:L_epsilon-made-explicit convex verzio}, which is surprisingly simple.
\begin{lem}
		\begin{equation*}
				L_\eps(\Om)
			\subset
				\begin{cases}
					\Om,
				&
					\text{if }\eps \geq \In(\Om),
				\\
					\Om
				\setminus
					\Int\left[
						\left(
							1-
							 \frac
							 	{\eps}
							 	{\In(\Om)}
						\right)
					\cdot
						\Om
					\right],
				&
					\text{if }\eps < \In(\Om);
				\end{cases}
		\end{equation*}
	where $x \mapsto \la\cdot x$ denotes enlargement, scale factor $\la$, about a fixed incentre of $\Om$.
\end{lem}

\begin{proof}
	We first show $f : x \mapsto d(x,\partial\Om)$ is concave on $\Om$, as follows. Suppose that $B_0$, $B_1$ are the largest balls in $\Om$ based respectively at $x_0$, $x_1 \in \Om$ (i.e., with radii $f(x_0)$ and $f(x_1)$). Then, for any $\la\in(0,1)$,
		$$
			B_\la := \la B_1 + (1-\la) B_0
		$$
	is a ball based at $x_\la := \la x_1 + (1-\la) x_0$, with radius $\la f(x_1) + (1-\la)f(x_0)$, contained wholly in $\Om$ by convexity. Hence,
		$$
				f(x_\la)
			\geq
				\la f(x_1) + (1-\la)f(x_0),
		$$
	as claimed.
	Therefore, given $y \in \Int(\Om)$ and $\eps < f(y)$, contracting $x\in \Om$ towards $y$ by a factor $\left(1 - \frac\eps{f(y)}\right)$ takes it outside of the $\eps$-inner neighbourhood:
		$$
				f 
				\left(
					\left(
						\frac
							\eps
							{f(y)}
					\right)
						y
					+
					\left(
						1
					-
						\frac
							\eps
							{f(y)}
					\right)
					x
				\right)
			\geq
					\left(
						\frac
							\eps
							{f(y)}
					\right)
					f(y)
				+
					\left(
						1
					-
						\frac
							\eps
							{f(y)}
					\right)
						f(x)
			\geq \eps.
		$$
	%
	Thus taking $y$ such that $f(y) = \In(\Om)$ gives the second case of \eqref{eq:L_epsilon-made-explicit convex verzio}. (The first is trivial.)
\end{proof}

\begin{figure}[ht]

	\centering
	\begin{tikzpicture}[scale = 6]

		\path[fill = blue!5] (-0.00616, 0.08352) rectangle (0.85902,0.91987);

		\path[fill = blue!5] (-0.00616, 0.08352) ++ (1,0) rectangle (1.85902,0.91987);

		\draw[semithick,fill = white] (0.04384, 0.77957) --
				(0.10215, 0.16202) --
				(0.39443, 0.13352) --
				(0.69521, 0.36598) --
				(0.80902, 0.72602) --
				(0.79672, 0.76476) --
				(0.29577, 0.86987) --
				(0.04384, 0.77957) --
				cycle;

		\node (B) at (0.386, 0.527) [circle,fill=black,inner sep = 0.5] {};

		\draw[fill = yellow, opacity = 0.1] (B) circle[radius = 0.317];

		\draw[semithick, fill = white] (1.04384, 0.77957) --
				(1.10215, 0.16202) --
				(1.39443, 0.13352) --
				(1.69521, 0.36598) --
				(1.80902, 0.72602) --
				(1.79672, 0.76476) --
				(1.29577, 0.86987) --
				(1.04384, 0.77957) --
				cycle;

		\node (A) at (1.386, 0.527) [circle,fill=black,inner sep = 0.5] {};

		\draw[fill = yellow, opacity = 0.1] (A) circle[radius = 0.317];

		\foreach \x in {{(0.04384, 0.77957)},{(0.10215, 0.16202)},{(0.39443, 0.13352)},{(0.69521, 0.36598)},{(0.80902, 0.72602)},{(0.79672, 0.76476)},{(0.29577, 0.86987)},{(0.04384, 0.77957)}}
		{
			\draw \x ++ (1,0) -- (A);
		}

	 	\draw (0.10215, 0.16202) -- (0.324, 0.383) -- (B) -- (0.324, 0.383) -- (0.39443, 0.13352);

	 	\draw (0.80902, 0.72602) -- (0.762, 0.722) -- (0.443, 0.54) -- (B) -- (0.443, 0.54) -- (0.69521, 0.36598) -- (0.443, 0.54) -- (0.762, 0.722) -- (0.79672, 0.76476);

	 	\draw (0.04384, 0.77957) -- (0.314, 0.607) -- (0.29577, 0.86987) -- (0.314, 0.607) -- (B);

		\path (-0.00616, 0.08352) -- (0.85902,0.08352) node[midway, anchor = north, inner sep = 5pt] {$\{\Om_S\}$};

		\path (-0.00616, 0.08352) ++ (1,0) -- (1.85902,0.08352) node[midway, anchor = north, inner sep = 5pt] {$\{\De_S\}$};

	\end{tikzpicture}
	\caption{The decompositions of convex polygon $\Om$ into $\{\Om_S\}$ and $\{\De_S\}$. The largest inscribed ball is depicted in pastel yellow.}
	\label{fig:1}
\end{figure}

We finally prove the generalisation of Heron's formula, as below. Both bounds are in fact tight.
\begin{prop}
	%
		$$
				\frac
					{\vol(\Om)}
					{\per(\Om)}
			\le
				\In(\Om)
			\le
				n
				\frac
					{\vol(\Om)}
					{\per(\Om)}.
		$$
	%
\end{prop}

Since any convex body can be approximated arbitrarily well by convex polytopes, we will prove the result assuming $\Om$ is a polytope. The proof then relies on decomposing $\Om$ in two different ways, illustrated in Figure \ref{fig:1}, which coincide when $\Om$ is a simplex.

\begin{proof}[Proof of lower bound]
		Consider one side $S\subset \partial\Om$ of $\Om$, and recall the locus $\Om_S$ of points in $\Om$ whose closest side of $\Om$ is $S$. Treating $S \subset \RR^{n-1}\times \{0\}$ and $\Om \subset \RR^{n-1}\times [0,\infty)$, we see that
			\begin{itemize}

				\item every point in $\Om_S$ can be at most $\In(\Om)$ distance from $S$, and

				\item the locus is bound by the hyperplanes bisecting $S$ and each of the neighbouring sides:

			\end{itemize}
		by convexity, the interior angle between any two sides lies in $(0,\pi)$, and thus any bisector with $S$ makes an interior angle at most $\pi/2$ with $S$: i.e., is either vertical or sloping inwards towards $S$.		Hence  $\Om_S \subset S \times [0,\In(\Om)]$ and thus
			$$
					\vol(\Om_S)
				\leq
					\vol_{n-1}(S)
					\cdot
					\In(\Om).
			$$
		Summing over $S$ gives the required inequality.
\end{proof}
\begin{proof}[Proof of upper bound]
		Choose an incentre $x\in \Om$.
		Fix side $S$, and let $\De_S$ denote the cone with base $S$ and vertex $x$. Again take $S \subset \RR^{n-1}\times \{0\}$ and $\Om \subset \RR^{n-1}\times [0,\infty)$ for simplicity.
		Since the ball of radius $\In(\Om)$ based at $x$ is contained in $\Om \subset \RR^{n-1}\times [0,\infty)$,
		we must have $x \in \RR^{n-1}\times [\In(\Om),\infty)$.
		Therefore, since $\De_S$ is a cone with base $S$ and height at least $\In(\Om)$,
			$$
					\frac1n
					\vol_{n-1}(S)\cdot\In(\Om)
				\leq
					\vol(\De_S)
				.
			$$
		The proof is finished by summing over $S$.
\end{proof}

\end{document}